\documentclass[11pt]{amsart}

\newtheorem{thm}{Theorem}[section]
\newtheorem{prop}[thm]{Proposition}
\newtheorem{lem}[thm]{Lemma}
\newtheorem{cor}[thm]{Corollary}

\theoremstyle{definition}

\newtheorem{remark}[thm]{Remark}
\newtheorem{example}[thm]{Example}

\newtheorem{algorithm}[thm]{Algorithm}

\newcommand{\bZ}{{\mathbb{Z}}}

\newcommand{\bP}{{\mathbb{P}}}
\newcommand{\bA}{{\mathbb{A}}}

\newcommand{\ine}{{\bf in}}

\newcommand{\Char}{\mbox{\rm char\,}} %% \char is already a command

\newcommand{\Gal}{\operatorname{Gal}}

\newcommand{\lra}{\longrightarrow}

\newcommand{\Sym}{{\operatorname{S}}}
\newcommand{\T}{{\operatorname{T}}}
\newcommand{\tr}{\mbox{\rm tr}}
\newcommand{\Ri}{k[X_1^{\pm 1}, \dots, X_m^{\pm 1}]_0}

\begin{document}

\title[symmetric functions]{Symmetric functions
and the phase problem in crystallography}
\author[J. Buhler and Z. Reichstein]{J. Buhler and Z. Reichstein}
\address{Reed College,
Portland, OR 97202}
\email{jpb@@reed.edu}
\address{Department of Mathematics, University of British Columbia,
Vancouver, BC V6T 1Z2, Canada}
\thanks{Z. Reichstein was partially supported by NSF grant DMS-901675
and by an NSERC research grant}
\email{reichst@@math.ubc.ca}

%%%%%%%%%%%%%%%%%%%%%%%%%%%%%%%%%%%%%%%%%%%%%%%%%%
% 05E05 Symmetric functions
% 13A50 Group actions on rings, invariant theory
% 20C10 Integral representations of finite groups
% 13P99 Computational aspects of commutative algebra
%%%%%%%%%%%%%%%%%%%%%%%%%%%%%%%%%%%%%%%%%%%%%%%%%%

\begin{abstract}
The calculation of crystal structure from X-ray diffraction data
requires that the phases of the ``structure factors'' (Fourier
coefficients) determined by scattering be deduced from the absolute
values of those structure factors.
Motivated by a question of Herbert Hauptman,
we consider the problem of determining phases by direct algebraic
means in the case of crystal structures with $n$ equal atoms in the
unit cell, with $n$ small.  We rephrase the problem as a question
about multiplicative invariants for a particular finite group action.
We show that the absolute values form a generating set for the field
of invariants of this action, and consider the problem of making this
theorem constructive and practical; the most promising approach for
deriving explicit formulas uses SAGBI bases.  \end{abstract}

\maketitle
% \tableofcontents

\section{Introduction}
\label{sect.intro}

If the unit cell of a crystal has $n$ atoms, located
at positions ${\bf r}_j$, $1 \le j \le n$, then the structure factor
associated to a reciprocal lattice vector ${\bf v}$ is
$$E_{\bf v} := \sum_{j=1}^n a_j \exp(2\pi i {\bf v} \cdot {\bf r}_j) \, , $$
where the $a_j$ are the scattering amplitudes determined by the
electron charge distribution in the $j$-th atom.
This is, in effect, a Fourier transform coefficient, and
the structure of the crystal can be determined from the $E_{\bf v}$
by an inverse Fourier transform.
However, in standard diffraction experiments, it is impossible to
measure the
$E_{\bf v}$ -- only their absolute values are observable.
The ``phase problem'' of crystallography
is to determine the phases of the $E_{\bf v}$ given magnitudes
$|E_{\bf v'}|$ for sufficiently many~{$\bf v'$}; this problem is
fundamental in the subject, and has received considerable
attention (\cite{giacov}, \cite{haupt.phase}).  The problem
of retrieving phase information from absolute values, together
with other physical constraints, occurs in several other areas
of physics, astronomy, and engineering.

In the crystallographic context, the phases can be determined in principle
only up to an additive constant, which is equivalent, via a Fourier transform,
to the indeterminacy of the origin of the crystal.  It is natural to
consider ``structure invariants,''  which are multiplicative combinations
of structure factors that are invariant under change of origin, i.e.,
additive translation of the phase.  In addition, structure invariants
play an important role in commonly used stochastic methods for phase retrieval.
The function
\begin{equation} \label{e.phase-m}
E_{\bf v_1} E_{\bf v_2} \dots E_{\bf v_{m}}
\end{equation}
is easily seen to be a structure invariant when
the reciprocal lattice vectors ${\bf v_i}$ sum to zero.
The most common case is $m = 3$, i.e., the
``triplet-structure invariant''
$$E_{\bf v_1} E_{\bf v_2} E_{-{\bf v_1} - {\bf v_2}}.$$

Exact formulas for phases of triplet structure invariants
are known in terms of magnitudes
for $n = 1,2,3$ \cite{haupt.algdirect}.
Herbert Hauptman asked one of us for a formula for
arbitrary~$n$, and the purpose of this paper is to
show that, at least in the case in which all atoms in the crystal
have the same $a_j$, such a formula exists,
and to explore techniques for finding such formulas explicitly.

In this paper we shall only consider crystals with
equal atoms (or equal polyatomic clumps); we will set the identical
scattering factors $a_j$ equal to~1.  In addition,
it is convenient to assume that the space group is the most basic
group P1 (isomorphic to the group $\bZ^3$ of translations);
see \cite[Appendix 1]{haupt.algdirect} for a description
of how to generalize to arbitrary space groups.

Under these assumptions the triplet phase determination problem can
be converted to a question about multiplicative invariant functions.
Somewhat to our surprise, this question seems to be new.

To express the triplet phase problem to a problem in symmetric functions,
start by noting that
the phase $\phi_{{\bf v_1},{\bf v_2}}$ of the triplet structure factor
$E_{\bf v_1} E_{\bf v_2} E_{-{\bf v_1} -{\bf v_2}}$ satisfies
\[ E_{\bf v_1} E_{\bf v_2} E_{-{\bf v_1} - {\bf v_2} } =
\left|E_{\bf v_1}\right| \, \left|E_{\bf v_2}\right| \,
\left|E_{- {\bf v_2} - {\bf v_2}}\right| \, \exp(i\phi_{{\bf v_1},{\bf v_2}})\,.
\]
Thus the cosine of the phase can be expressed in terms of absolute
values and the sum of the triplet-structure invariant and its
complex conjugate.

We want to express the phase $\phi_{{\bf v_1},{\bf v_2}}$ in terms
only of absolute values
$$\left| E_{a{\bf v_1} + b{\bf v_2}} \right|^2 \, = \,
E_{a{\bf v_1} + b{\bf v_2}} E_{a{\bf v_1} + b{\bf v_2}}^*
\, = \, E_{a{\bf v_1} + b{\bf v_2}} E_{-a{\bf v_1} - b{\bf v_2}}$$
corresponding to reciprocal lattice vectors $a{\bf v_1} + b{\bf v_2}$,
where $a$ and~$b$ are integers.

Fix ${\bf v_1}$ and ${\bf v_2}$ and let
\begin{eqnarray*}
x_j = \exp(2\pi i {\bf v_1} \cdot {\bf r}_j), & \quad &
y_j = \exp(2\pi i {\bf v_2} \cdot {\bf r}_j), \quad 1 \le j \le n \\
X = (x_1 , \cdots , x_n), & \quad & Y = (y_1, \cdots , y_n) \; .
\end{eqnarray*}
Then
$$ E_{\bf v_1} = \sum_{j=1}^n x_j \quad \mbox{and} \quad
E_{\bf v_2} = \sum_{j=1}^n y_j \; .  $$
Using the fact that
$z+z^* = 2|z|\cos(\phi)$
if $z$ is a complex number with absolute value $|z|$ and argument~$\phi$,
we see that we need to express
\begin{eqnarray*}
E_2 (X, Y) &:= &
E_{\bf v_1} E_{\bf v_2} E_{- {\bf v_1} - {\bf v_2}}
 + E_{- {\bf v_1}} E_{- {\bf v_2}} E_{{\bf v_1} +{\bf v_2}} \\
\label{e.E_2}
 & = &
\left|E_{\bf v_1}\right| \; \left|E_{\bf v_2}\right| \;
\left|E_{- {\bf v} - {\bf v_2}}\right| \;
2 \cos(\phi_{{\bf v_1},{\bf v_2}}) \\
 & = & \sum_{j=1}^n x_j \sum_{k=1}^n y_k \sum_{l=1}^n \frac{1}{x_l y_l} \; + \;
\sum_{j=1}^n \frac{1}{x_j} \sum_{k=1}^n \frac{1}{y_k} \sum_{l=1}^n x_l y_l \\
 & = & \sum_{j,k,l=1}^n \frac{x_j y_k}{ x_l y_l} \; + \;
\sum_{j,k,l=1}^n \frac{x_l y_l} {x_j y_ k}\\
\end{eqnarray*}
in terms of the magnitudes
\begin{eqnarray*}
 q_{a , b}(x_1, \dots, x_n, y_1, \dots, y_n) & := &
E_{\bf v} E_{\bf -v} \\
& = & \sum_{j=1}^n x_j^a y_j^b
\; \cdot \; \sum_{j=1}^n x_j^{-a} y_j^{-b} = \sum_{i,j = 1}^n
\frac{x_i^a y_i^b}{x_j^a y_j^b} \, ,
\end{eqnarray*}
for suitable integers $a$ and $b$; here
${\bf v} = a {\bf v_1} + b {\bf v_2}$ is an arbitrary reciprocal lattice
vector.

We will call $q_{a , b}$ an ``observable'' since it is (the
square of) an absolute value and therefore it is possible to
observe it physically.  Thus our goal is to express $E_2$ as a
rational function of the observables $q_{a, b}$.
{}From now on we will treat this as a question about variables $x_i$
and $y_i$, ignoring the fact that they are complex numbers of
absolute value one. A simple Zariski density argument shows that
this does not change the underlying problem, i.e., we are not going to
``miss" any identities by assuming that $x_i$ and $y_i$ are arbitrary
complex numbers, rather than just those of absolute value one.

\begin{example}
The reader can easily verify that for $n = 2$ the triplet phase invariant
is a polynomial in three observables:
$$E_2(X,Y)  = 2(q_{1,0}+q_{0,1}+q_{1,1})-8.$$
\end{example}
\begin{example} \label{first.threeformula}
The formula for $n = 3$ is considerably more elaborate, and was
discovered by Hauptman \cite{haupt.algdirect}; it takes the form
$E_2(X,Y) = N/D$, where
$$D := q_{0,1}+q_{1,0}+q_{1,1}-3 = -3 + \sum_{i,j=1}^3
\left(\frac{x_i}{x_j} + \frac{y_i}{y_j} +
\frac{x_i y_i}{x_j y_j}\right)$$
and
\begin{eqnarray*}
N & := & 135-31D+D^2+2(q_{1,0}q_{0,1}+q_{1,0}q_{1,1}+q_{0,1}q_{1,1})\\
&&+(q_{0,1} q_{2,1}+ q_{1,0} q_{1,2}+ q_{1,-1} q_{1,1})\\
&&-5 (q_{1,2} + q_{2,1}+ q_{1,-1}) -2( q_{0,2}+ q_{2,2} + q_{2,0}) \, .
\end{eqnarray*}
(For another formula, see Example~\ref{ex.formula3}.)
\end{example}

Higher order generalizations of triplet structure invariants
are defined as a product of structure factors $E_v$ where the
$v$ sum to~0.  If a structure invariant is the product of $m$
factors~\eqref{e.phase-m} for arbitrary~$m$, then this gives
rise, in a similar manner, to the problem of expressing
\begin{equation} \label{e.E_m}
E_m(X_1, \cdots , X_m) :=
\sum_{j, j_1, \dots, j_m = 1}^n \Bigl( \;
\frac{x_{1j_1} x_{2j_2} \dots x_{mj_m}}{x_{1j} \dots x_{mj}}  +
\frac{x_{1j} \dots x_{mj}}{x_{1j_1}x_{2j_2} \dots x_{mj_m}} \Bigl)
\end{equation}
as a rational function in observables
\begin{equation} \label{e.observables}
q_{r_1, \dots, r_m} := \sum_{i, j=1}^n
\frac{x_{1i}^{r_1} \dots x_{mi}^{r_m}}{x_{1j}^{r_1} \dots x_{mj}^{r_m}}
\, , \end{equation}
where $r_1, \dots, r_m$ are integers.

Note that the rational function $f = E_m$ has
the following invariance properties.
\begin{description}
\item[a] $f$ is of weight 0 in each $n$-variable vector
\[ X_1 = (x_{11}, \dots, x_{1n}), \quad \dots, \quad
X_m = (x_{m1}, \dots, x_{mn}) \, . \]
That is,
$f(c_1X_1, \dots, c_m X_m) = f(X_1, \dots, X_m)$
for any non-zero scalars $c_1, \dots, c_m$, or,
more succinctly: $f(c_jx_{ij}) = f(x_{ij})$.
\item[b] $f$ is self-reciprocal in the sense that
it remains unchanged if every variable $x_{ij}$ is
simultaneously replaced by $x_{ij}^{-1}$, i.e.,
$f(x_{ij}^{-1}) = f(x_{ij})$.
\item[c]
$f$ is multi-symmetric in the sense that it remains unchanged
if the variables in each array $X_i$ are (simultaneously)
permuted by the same permutation $\sigma \in \Sym_n$, i.e.,
$f(x_{i\sigma(j)}) = f(x_{ij})$.
\end{description}

Our main results are summarized in the following theorem.
Both parts answer questions posed by H. Hauptman.

\begin{thm} \label{thm0.1}
(1) (Theorem~\ref{thm2} or Theorem~\ref{thm3})
Every rational function $f(x_{ij})$ in $mn$ variables,
satisfying the invariance properties
{\bf a}, {\bf b}, {\bf c}, can be expressed as a rational function
of the observables $q_{c_1, \dots, c_m}$.
In particular, $E_m$ can be expressed as a rational function
of the observables $q_{c_1, \dots, c_m}$ for any $m \geq 1$.

\smallskip
(2) (Proposition~\ref{prop.n>=5}(b))
Suppose $n \geq 4$. Then
$E_m$ is not a polynomial in the observables
$q_{c_1, \dots, c_m}$ for any $m \geq 2$.
\end{thm}
Our proof of Theorem~\ref{thm0.1}
requires the consideration of all $f$ satisfying
the invariance properties {\bf a}, {\bf b}, {\bf c}, even if one
is only interested in the structure invariants~$E_m$.
This is true of most of our other proofs and algorithms.

The overall outline of the paper is as follows.
In Section~\ref{sect1} we prove (a slightly more precise
version of) Theorem~\ref{thm0.1}(2)
for $m = 1$ using basic Galois theory and some combinatorial
arguments; the proof is not obviously constructive.
In Section~\ref{sect.alg-single} we consider two possible
approaches to making the argument constructive.
In the subsequent section we give a fast algorithm for $n \leq 4$,
based on SAGBI bases.
In Section~\ref{sect.mult} we turn our attention
to the multi-array case, i.e., $m>1$.
In Section~\ref{sect.explicit-phase}
we reduce the problem of computing the invariant $E_m$ to that of
expressing certain invariant functions in terms of observables
in the single array case ($m = 1$).  We then use the SAGBI
basis algorithm of Section~\ref{sect.sagbi}
to obtain new expressions for $E_2$ in the case where $n = 3$
and $4$; see Examples~\ref{ex.formula3} and~\ref{ex.formula4}.
In Section~\ref{sect.regular}
we prove Theorem~\ref{thm0.1}(2).
Finally, in Section~\ref{sect.rationality}
we study the structure of the field of rational functions $f$
satisfying the invariance properties ({\bf a}) - ({\bf c}) as
an abstract field, without reference to the observables.

In the course of our work on this paper, we have encountered
a phenomenon that often arises in the
interstices between mathematics and its applications. Depending
on the context, solving a mathematical problem can mean many
different things, e.g.,

\smallskip
(a) proving a theorem,

\smallskip
(b) giving a constructive proof, or

\smallskip
(c) giving an algorithmic proof, suitable for practical computations.

\smallskip
\noindent
As one moves down this list, the problem can become
more difficult, requiring different
techniques and ways of thinking.  However, there is
usually a subtle but important interplay between these
different modes of solution.

\section*{Acknowledgments} The authors would like to thank
H.~Hauptman for bringing this problem to their attention,
M.~Lorenz for contributing Proposition~\ref{prop.rationality}, and
D.~Eisenbud, N.~Elkies, J.~Friedman, I.~Laba, H.~Lenstra, G.~Martin,
D.~Peterson, C.~Procesi, I.~Swanson, and R.~Thomas for helpful conversations.

\section{One set of variables}
\label{sect1}

Fix a field $k$ of characteristic~0.
Let $X = (x_1, \dots, x_n)$ be an $n$-tuple of independent variables
over $k$.  We will operate on $n$-tuples as if they were diagonal matrices, so
$X^{-1} = (x_1^{-1}, \dots, x_n^{-1})$, $\tr(X) = x_1 + \dots + x_n$, etc.

Let $k(X)_0 \subset k(x_1, \cdots, x_n)$
be the field of rational functions in the $x_i$ of total degree 0; in
other words, an element $f\in k(X)_0$ is a quotient of
homogeneous polynomials of the same degree.
Equivalently, $k(X)_0$ is the field generated by the $x_i/x_j$:
$$ k(X)_0 = k(x_i/x_j : 1 \le i, j \le n, i \ne j )\; .$$
We note that the field $k(X)_0$ can also be viewed as
the function field of the projective space $\bP^{n-1}$.

The symmetric group $\Sym_n$
acts on $k(X)_0 = k(x_1, \cdots , x_n)_0$ by permuting the variables
$x_1, \cdots, x_n$ in
the natural way.  In addition, we let $\tau$ denote the
automorphism that takes $X$ to $X^{-1}$, i.e.,
$$ \label{e.tau}
\tau (x_i) = \frac{1}{x_i}, \quad 1 \le i \le n.
$$
This automorphism is obviously of order two, and we let $\T = \{1,\tau\}$
denote the corresponding group.  The actions of $\Sym_n$ and $\T$ commute
so that the group
$$ G := \Sym_n \times \T $$
acts on $k(X)_0$.
This action is faithful for $n \geq 3$.
If $n = 2$, $G = \Sym_2 \times \T$ has order 4, and
the kernel of its action on $k(x_1, x_2)_0$ is the subgroup
of order 2 generated by $(\sigma, \tau)$, where $\sigma$ is the nontrivial
element of $\Sym_2$.

The main theorem of this section is that the observables
\begin{equation} \label{e.q_r}
q_r := \tr(X^r) \tr(X^{-r}) =
\sum_{i,j = 1}^n \left(\frac{x_i}{x_j}\right)^r
\end{equation}
generate the invariant field $k(X)_0^G$. (Note that here $m = 1$, so
that the observables~\eqref{e.observables} have only one subscript.)

\begin{thm} \label{thm1}
$$k(X)_0^G =
k(q_r \, | \, 1 \le r \le n(n-1)/2) \; .
$$
\end{thm}

\smallskip
Before proving the theorem, we prove two lemmas, one combinatorial and
the other algebraic.

Let $n \geq 2$ be an integer, $N = n(n-1)$, and $\Lambda$ be the $N$-element
set
\begin{equation} \label{e.Lambda}
\text{$\Lambda = \{ (i, j)$ $|$ $i, j = 1, \dots, n$ and $i \neq j \}$.}
\end{equation}
When convenient we will sometimes omit the comma
and write $(ij)$ instead of $(i,j)$.
In addition, we tend to visualize the elements of $\Lambda$ as the
off-diagonal positions in an $n$ by~$n$ matrix.

There is a natural action of $G = \Sym_n \times \T$ on $\Lambda$, where
$\tau$ acts by transposition
$$\tau(i,j) = (j,i)$$
and the symmetric group $\Sym_n$ acts simultaneously on the rows and
columns:
$$\sigma (i,j) = (\sigma(i),\sigma(j)), \quad \sigma \in \Sym_n.$$
This gives a map from $G$ to the symmetric group
$\Sym_N = \Sym_{n^2-n} = \text{Sym}(\Lambda) $ of
all permutations of~$\Lambda$, and this map is an injection for
$n \ge 3$, in which case we will usually just choose to regard $G$ as
a subgroup of~$\Sym_{N}$.  In the case $n = 2$ the map is
surjective with kernel of order two as described earlier.

We will say that elements $x = (i, j)$ and $x' = (i', j')$ of $\Lambda$
are {\em opposite} if $x$ is the transpose of $x'$, i.e.,
$i = j'$ and $j = i'$.  If exactly one of these equalities holds,
i.e., $x'$ is not opposite to~$x$, but it lies in the same row or column
as the transpose of $x$, then we say that $x$ and $x'$ are {\em adjacent}.
Note that opposite pairs are not also adjacent.

We will say that $g \in \Sym_{N}$ preserves adjacency
(respectively opposition) if for any adjacent (respectively, opposite)
elements $x, x' \in \Lambda$, the images
$g(x)$ and $g(x')$ are also adjacent (respectively, opposite).
Our combinatorial lemma says that for $n > 2$, $G = \Sym_n \times \T$
is precisely the subgroup of $\Sym_N$ that preserves both of these relations.

We note that for $n = 2$ the situation is simple: $G$ maps
onto $\Sym_N = \Sym_2$ and the nontrivial element
of $\Sym_N$ preserves both adjacency and opposition.
So from now on we consider $n \ge 3$.

\begin{lem} \label{lem.adj}
Let $n \ge 3$.  Then
$h  \in \Sym_{N}$
preserves both adjacency and opposition if and only if
$h \in G \subset \Sym_N$.
\end{lem}

\begin{proof}
It is immediate from the definition that every
$h $ in~$G$ preserves both adjacency and opposition, so we
only need to prove that any element preserving these relations
lies in~$G$.

Suppose $h \in \Sym_{N}$
preserves both adjacency and opposition. To show that $h \in G$,
we will multiply $h$ by elements of $G$ until we arrive
at the identity permutation of $\Lambda$.

It is easy to see that $\Sym_n$ acts transitively on $\Lambda$.
Thus, after composing $h$ with an element
of $\Sym_n \subset G$, we may assume
$h(12) = (12)$.

We claim that we may also assume that $h(13) = (13)$.
Indeed, since $h$ preserves opposition, $h(12) = (12)$ implies
$h(21) = (21)$.  Suppose $h(13) = (ij)$. Since $(21)$ and $(13)$ are adjacent,
so are $(21)$ and $(ij)$, i.e., either $i = 1$ or $j = 2$. If $i = 1$
then $j \geq 3$; thus after replacing $h$ by $[3, j] h$,
we obtain $h(12) = (12)$ and $h(13) = (13)$, as desired.
(Here $[3, j]$ denotes
the transposition in $\Sym_n$ that interchanges $3$ and $j$.)
On the other hand, if $j = 2$ then $i \neq 1, 2$ and, after
replacing $h$ by $[3, i] [1, 2] \tau h$, we once again obtain
$h(12) = (12)$ and $h(13) = (13)$.
This proves the claim.

Since $h$ preserves opposition, $h(13) = (13)$
implies that $h(31) = (31)$.  Now, since $(23)$ is the
unique pair adjacent to both $(31)$ and $(12)$, and $h(23)$ is
the unique pair adjacent to both $h(31) = (31)$ and $h(12) = (12)$,
we conclude that $H(23) = (23)$.  Since $h$ preserves opposition,
we also have $h(32) = (32)$.  Thus $h$ fixes $(ij)$ for $1 \le i,j \le 3$.
This completes the proof of Lemma~\ref{lem.adj} for $n = 3$; from now on we
will assume that $n \geq 4$.

Suppose $h(1i) = (ab)$ for some
$i \geq 4$. Since $(21)$ and $(1i)$ are adjacent, so are
$(21)$ and $(ab)$. That is, either
$a = 1$ or $b = 2$. Repeating this
argument with $(31)$ in place of $(21)$, we see that either
$a =1$ or $b = 3$. Since $b$ cannot be equal to both $2$ and $3$,
we conclude that $a = 1$. In other words, $h(1i) = (1 {\sigma(i)})$,
where $\sigma$ is a permutation of $4, 5, \dots, n$. After replacing
$h$ by $\sigma^{-1} h$, we reduce to the case where $h$ fixes $(1i)$
(and thus $(i1)$), for every $i = 2, \dots, n$.

We claim that $h$ is the identity permutation, i.e., that
$h(ab) = (ab)$ for every
$(ab) \in \Lambda$. Since we know this in the cases where $a = 1$ or $b = 1$,
we may assume $a, b \geq 2$. In this case $(ab)$ is the unique element of
$\Lambda$ that is adjacent to both $(1a)$ and $(b1)$. Hence, $h(ab) = (ab)$,
as claimed.
\end{proof}

Next, we prove an algebraic lemma from which the theorem will
follow easily.  Let
$$f(t) = \prod_{\begin{array}{c} \text{\tiny{$i,j = 1$}} \\
\text{\tiny{$i \neq j$}}
\end{array}}^n
\left (t - \frac{x_i}{x_j} \right ) =
\sum_{i=0}^{N} (-1)^i c_i t^{N-i}$$
be the polynomial of degree~$N =n(n-1)$ whose roots are the $x_i/x_j$,
$i \ne j$.
The coefficients $c_i$ are the elementary symmetric functions in
those roots, and since the reciprocals of roots are themselves
roots, the polynomial $f$ satisfies $t^Nf(1/t) = f(t)$, which is
equivalent to
$$c_i = c_{N-i}, \quad 0 \le i \le N.$$

The elements $q_r = \tr(X^r) \tr(X^{-r})$ are symmetric
functions of the $x_i/x_j$ are hence polynomials in the
$c_i$; we let
$$K := k(q_r \, | \, r = 1, 2, \cdots )$$
be the field generated by the observables~$q_r$.

As we will see, the proof of Theorem~\ref{thm1}
basically comes down to determining the Galois group of~$f$ over
the field~$K$.

\begin{lem} \label{lem2.1}
With the above notation,

\smallskip
(a) $k[c_1, \dots, c_r] = k[q_1, \dots q_r]$
for any $r$,  $ 1 \le r \le N$.

\smallskip
(b) $K = k(c_1, \dots, c_{N/2}) = k( q_1 , \dots , q_{N/2}).$

\smallskip
(c) $k(X)_0$ is the splitting field of $f(t)$ over $K$.
\end{lem}

\begin{proof}
Since
$$q_r - n = \tr(X^r) \tr(X^{-r}) - n =
\sum_{i \neq j} (x_i/x_j)^r$$
is the sum of the $m$-th powers of the
roots of $f(t)$, part (a) follows from Newton's formulas that
express the symmetric polynomials $c_1, \cdots , c_r$ in terms
of the power sums $q_1, \cdots , q_r$.

Part (b) follows from part (a) and the symmetry of the $c_i$.
Part (c) follows from (b) and the fact that $k(X)_0 = k(x_i/x_j)$.
\end{proof}

We are now ready to finish the proof of Theorem~\ref{thm1}.
Clearly, $K \subset k(X)_0^G$. Consider the tower
\[ \begin{array}{c}   k(X)_0 \\
                          |   \\
                     k(X)_0^G \\
                         |   \\
                         K
\end{array} \]
of field extensions. By the lemma, $k(X)_0$ is a Galois extension
of $K$.  Identify the set of roots of $f(t)$
with the set $\Lambda = \{(i,j) \, | \, i \neq j \}$, letting
$x_i/x_j \leftrightarrow (i,j)$.
The action of $G = \Sym_n \times \T$ on the set of roots
is the same as its action on $\Lambda$ described earlier.

The notions of adjacency and opposition in $\Lambda$
have a natural interpretation in this context.
If $x = (ab)$ and $x' = (cd)$ are elements of $ \Lambda$ let
$r = x_a/x_b$ and $r' = x_c/x_d$ be the
corresponding roots of $f(t)$. Then $x$ and $x'$ are adjacent if and only if
$r r'$ is again a root of $f(t)$ and opposite if and only if $r r' = 1$.
Thus any $\Gal(k(X)_0/K)$ acts on the set of roots in a way that preserves both
adjacency and opposition.  Lemma~\ref{lem.adj} now tells us
that $\Gal(k(X)_0/K) = G$.  Thus $k(X)_0^G = K$. This completes the proof of
Theorem~\ref{thm1}.

\begin{remark} \label{rem2.2}
Since $\Sym_n \subset G$ acts transitively on the roots
of $f(t)$, we conclude that $f(t)$ is irreducible over $k(X)_0^G$.
\end{remark}

\section{Constructive proofs}
\label{sect.alg-single}

The proof in the last section ultimately relies on a fundamental
and beautiful result in Galois theory: anything fixed by all elements
of a Galois group lies in the ground field. We will now
discuss a constructive proof, which could be viewed as the
result of tracing through the argument
of the previous section, rendering the underlying Galois
theory explicit at each step.  We note that both proofs
rely on Lemma~\ref{lem.adj}.

We begin by letting
$k[X^{\pm 1}]_0$ be the $k$-algebra whose elements are $k$-linear
combinations of Laurent monomials $x_1^{a_n} \dots x_n^{a_n}$
of total degree~0, i.e., where the $a_i$ are integers whose
sum is~0.  Note that $k[X^{\pm1}]_{0}$ is generated,
as a $k$-algebra, by elements of the form $x_i/x_j$; in particular,
the field $k(X)_0$ is the field of fractions of $k[X^{\pm 1}]_0$. Note also
that $k[X^{\pm 1}]_0$ is a $G$-invariant subring of $k(X)_0$.

Let $z_{ij}$ be a set of $N = n(n-1)$ algebraically independent
variables over $k$, where $i$ and $j$ are distinct integers between
$1$ and $n$.  For notational convenience,
we also set $z_{ii} = 1$ for every $i = 1, \dots, n$.
We now define a surjective $k$-algebra
homomorphism
\[ \phi \colon k[z_{ij}] \lra k[X^{\pm 1}]_0  \]
by $\phi(z_{ij}) = x_i/x_j$.
Let $s_r$ be the $r$th elementary
symmetric polynomial in the $N$ variables $z_{ij}$
and $p_r = \sum_{i \neq j} z_{ij}^r$ be the sum
of the $r$th powers of these variables.
Here $s_0 = 1$ and $\phi(s_r)$ is the element of $k[X^{\pm 1}]_0$
we called $c_r$ in the statement of Lemma~\ref{lem2.1}.
We let $\Sym_{N}$ denote the group of all permutations of
$\{ z_{ij}: i \ne j \}$, and identify $G = \Sym_n \times \T$
with the subgroup that acts on the $z_{ij}$ by
$\sigma(z_{ij}) = z_{\sigma(i) \sigma(j)}$, for $\sigma \in \Sym_n$, and
$\tau (z_{ij}) = z_{ji}$.

Define three polynomials $D_1, D_2, D \in k[z_{ij}]$ of $N$ variables
by
\[ D_1(z_{ij}) =
 \prod_{\substack{\text{$(a b)$ and $(cd) \in \Lambda$} \\
\text{are not opposite}}}
(z_{ab} z_{cd} - 1) \, , \]
\[ D_2(z_{ij}) =
 \prod_{\substack{\text{$(ab)$ and $(cd) \in \Lambda$} \\
\text{are not adjacent}}}
(z_{ab} z_{cd} - z_{ef}) \, .  \]
and $D = D_1 D_2$.
As we shall see below, $\phi(D)$ is a ``universal
denominator'', such that if $f \in  k[X^{\pm 1}]_0^G$
then $\phi(D)f$ is a polynomial in the observables.

\begin{lem} \label{lemD}
(a) Suppose $g \in G$. Then $\phi(gD) = \phi(D) \neq 0$

\smallskip
(b) Suppose $g \not \in G$. Then $\phi(gD) = 0$.
\end{lem}

\begin{proof} (a) Clearly $gD_1 = D_1$ and $gD_2 = D_2$,
because $g$ preserves both adjacency and opposition. Thus $gD = D$,
so that $\phi(gD) = \phi(D)$. To show that $\phi(D) \neq 0$, note that
if the image under $\phi$ of a factor of $D_1$ is zero then
$$\frac{x_a}{x_b} \frac{x_c}{x_d} = 1\,.$$
This implies that $(ab)$ and $(cd)$ are opposite elements of $\Lambda$,
which is excluded by the definition of $D_1$.
This shows that $\phi(D_1) \neq 0$.
Similarly, $\phi(z_{ab} z_{cd} - z_{ef}) = 0$ if and only if $(ab)$ and $(cd)$
are adjacent in $\Lambda$. Thus $\phi(D_2) \neq 0$, and consequently,
$\phi(D) \ne 0$.

\smallskip
(b) If $g \in \Sym_n$ is not in~$G$ then by Lemma~\ref{lem.adj}
at least one of the following holds: (i) $g$ does not preserve opposition
or (ii) $g$ does not preserve adjacency.

If (i) holds then $g^{-1}$ does not preserve opposition either.
In other words, there exists a pair of non-opposite elements $(ab)$
and $(cd)$ such that $g(ab)$ and $g(cd)$ are
opposite, say, $b_1  = c_1$. Then $z_{ab} z_{cd} -1$ is a factor of
$D_1$ and $\phi(g(z_{ab} z_{cd} - 1)) = 0$. Hence,
$\phi(g (D_1)) = 0$ and thus $\phi(g(D)) = 0$.

Similarly, if (ii) holds then
there exists a pair of non-adjacent elements $(ab)$ and
$(cd)$ such that $g(ab) = (a_1 b_1)$ and $g(cd) = (c_1 d_1)$ are adjacent,
say, $b_1 = c_1$. Now
\[ \phi(g(z_{ab} z_{cd} - z_{g^{-1}(a_1 d_1)})) =
\phi(z_{a_1b_1} z_{c_1 d_1} - z_{a_1 d_1})) = \frac{x_{a_1}}{x_{b_1}}
\frac{x_{c_1}}{x_{d_1}} - \frac{x_{a_1}}{x_{d_1}} = 0 \, \]
so that $\phi(g(D_2)) = 0$ and thus $\phi(g(D)) = \phi(g(D_1)) \phi(g(D_2))
= 0$, as claimed.
\end{proof}

With the universal denominator $\phi(D)$ in hand, we can now state
our algorithm.

\begin{algorithm} \label{alg.single}
\hfill

\begin{description}
\item[Input] A function $f(x_1, \dots, x_n)$ in $  k(X)_0^G$.
\item[Output] A rational function in the $q_r$ representing $f$.
\item[Step 1] Write $f = f_1/f_2$ of a quotient of elements
$f_i$ that are in $k[X^{\pm 1}]_0^G$;
apply each of the subsequent steps to
$f_1$ and $f_2$ (to simplify the notation we just
assume from now on that $f \in k[X^{\pm 1}]_0^G$).
\item[Step 2]  Find an element $F \in k[z_{ij}]$ such that $\phi(F) = f$.
\item[Step 3] Set $A = \sum_{g \in \Sym_{N}} g(DF)$, write the
numerator and denominator of
$$ \frac{\phi(A)}{2 n! \, \phi(D)} $$
as polynomials in the $q_r$ and output the result.
\end{description}
\end{algorithm}

\noindent We comment on each step in turn.
\smallskip

\noindent{\bf Step 1:}
To express an invariant degree 0 rational function $f$ as a quotient of
invariant polynomials, recall that
$f$ is a quotient of two polynomials
of equal degree, say of degree $d$. Dividing top and bottom by
$x_1^d$, we can write $f$ in the form $f_1 = a/b$, where
$a$ and $b$ are in $K[X^{\pm 1}]_0$ (but $a$ and $b$ may
not be $G$-invariant). Since $f$ is $G$-invariant,
$$f = \frac{1}{2n!} \sum_{g \in G} g(f) =
\frac{1}{2n!} \sum_{g \in G} \frac{g(a)}{g(b)} =
\frac{f_1}{f_2} $$
where the numerator and denominator $f_i$ are $G$-invariant polynomials
in $k[X^{\pm 1}]_0^G$.
It suffices to express $f_1$ and~$f_2$ as rational functions in the $q_r$, and
we can therefore assume that $f$ is in $k[X^{\pm 1}]_0^G$ from now on.

\smallskip
\noindent{\bf Step 2:}
To lift $f$ to an element $F \in k[z_{ij}]$,
write $f(x_1, \dots, x_n)$ as a $k$-linear combination
of degree-0 Laurent monomials
$x_1^{a_1} \dots x_n^{a_n}$, with $a_1 + \dots + a_n = 0$.
Any such monomial is a product of a finite number of terms of the form
$\frac{x_i}{x_j}$.  Write all of the monomials in~$f$ in this form and
replace each $\frac{x_i}{x_j}$ by $z_{ij}$ to
obtain the desired $F \in k[z_{ij}]$.

\smallskip
\noindent{\bf Step 3:}
By Lemma~\ref{lemD} the only non-zero terms in the sum
$$\phi(A) = \sum_{g \in \Sym_N} \phi(g(D))\phi(g(F))$$
correspond to $g \in G$. Thus
$$\phi(A) = \sum_{g \in G} \phi(g(D))\phi(g(F)) =
\sum_{g \in G} \phi(D)g(\phi(F)) = 2n! \, \phi(D) \, f \;,$$
as claimed.

The only remaining thing that needs to be done is to explain how
to write $\phi(A)$ and $\phi(D)$ as polynomials in the~$q_r$.
Clearly $A$ is a symmetric polynomial in the $z_{ij}$, where
$1 \le i, j \le n$.
It is therefore a polynomial in the elementary symmetric functions
$s_i$, $1 \le i \le N$, and a polynomial in the
power sums $p_r = \sum \, z_{ij}^r$ by Newton's formulas.
Since $\phi(p_r) = q_r-n$ it follows that $\phi(A)$ can be written
as a polynomial in the $q_r$, as desired.

{}From the formula
\[ 2 n! \phi(D)  =  \sum_{g \in G} \phi(g(D))
 =  \sum_{g \in \Sym_N} \phi(g(D)) =
\phi \left(\sum_{g \in \Sym_N} g(D)  \right) \]
it follows that, once $D$ is symmetrized over~$\Sym_N$, the same procedure
can be used to express $2 n ! \phi(D)$ in terms of the observables.

\begin{remark} \label{rem.alg1}
The above procedure can be modified to produce polynomials in the
$q_r$ with $r \le N/2$. Indeed, recall that
$\phi(s_{N-i}) = c_{N-i} = c_i = \phi(s_i)$.  Thus the image
under~$\phi$ of a polynomial $P(s_1, \cdots, s_N)$
in the elementary symmetric functions is unchanged if $s_i$ is
replaced by $s_{N-i}$ for $i > N/2$.

Note also that if one only wants to
express $f(x_1, \dots, x_n)$ as a rational function
in $q_r$ for $1 \le r \le N$, rather than $1 \le r \le N/2$, then
the algorithm of~\cite[Proof of Proposition 1.1.2]{sturmfels1}
can be used to write $A$ and the symmetrization of~$D$ directly as
polynomials in the power sums $p_i$
without going through elementary symmetric polynomials.
\end{remark}

\section{Initial terms and SAGBI bases}
\label{sect.sagbi}

It is natural to try to apply Gr\"obner bases to our problem.
Although Gr\"obner bases are usually applied
to rings, they can be adapted to solve problems in function fields.
This was first noticed in detail by Sweedler \cite{sweedler},
and has since been extended in several theoretical
and practical ways; see, e.g.,~\cite{muller}.
In our context we would have a $G$-invariant polynomial $f(X)$
in $x_i/x_j$ that, by Theorem~\ref{thm1}, lies
in the field $k(q_k) \subset k(X)_0$.
By introducing extra variables and calculating the Gr\"obner basis
of a suitably chosen ideal, an explicit expression can be found
for $f$ as a rational function in the~$q_i$.

Unfortunately, all of our implementations of this idea suggest
that it has the same trouble as implementations of the constructive
algorithm given in the previous section: they are far too slow.
% In both cases, they have trouble with
% the first nontrivial example \ref{first.threeformula} which can
% in fact be obtained (with some ingenuity) by hand computations.
The purpose of this section is to introduce a faster algorithm,
for $n \leq 4$, using a variant of Gr\"obner bases called
SAGBI bases.

We shall always assume that $f \in k[X]_0^G$;
the general case reduces to
this one (cf. Algorithm~\ref{alg.single}, Step 1).

\subsection*{The subduction algorithm}
Given an element
\[ p(x_1, \dots, x_n)  = \sum c_{a} x^{a} \in k[X^{\pm 1}] =
k[x_1^{\pm 1}, \dots, x_n^{\pm 1}] \, ,  \]
where $a = (a_1, \dots, a_n) \in \bZ^n$, $x^a = x_1^{a_1} \dots x_n^{a_n}$
and $c_a \in k$. We will write $\ine(p)$ for the initial exponent $p$, i.e.,
the lexicographically largest exponent $a$ such that $c_a \neq 0$.
If $R$ is a subalgebra of $k[X^{\pm 1}]$ then $\{ \ine(p) \, | \, p \in R \}$
is clearly a subsemigroup of $\bZ^n$; this semigroup is usually denoted by
$\ine(R)$. We are interested in the case
where $R = k[X^{\pm 1}]_0^G$; in this case
$\ine(R)$ consists of elements $a = (a_1, \dots, a_n) \in \bZ^n$
satisfying the following conditions:
\begin{equation} \label{e.S}
\begin{array}{l}
\text{(i) $a_1 + \dots + a_n = 0$,} \\
\text{(ii) $a_1 \geq \dots \geq a_n$, and} \\
\text{(iii) $(a_1, \dots, a_n) \succeq (- a_n, \dots, - a_1)$.}
\end{array}
\end{equation}
Here and in the sequel, $\succ$ denotes the lexicographic order on $\bZ^n$.

\begin{prop} \label{prop.subduct}
Suppose $B$ is subset of $R = k[X^{\pm 1}]_0^G$ chosen so that
the elements $\ine(b)$ generate $\ine(R)$ as a semigroup,
as $b$ ranges over $B$.  Then $R = k[B]$.
\end{prop}

Our proof below is based on the subduction algorithm of
Robbiano-Sweedler~\cite{rs} and Kapur-Madlener~\cite{km}
for expressing a given element
$\alpha \in R$ as a polynomial in elements of $B$.

\begin{proof} We want to write $\alpha \in R$ as a polynomial in
elements of $B$. If $\alpha = 0$, we are done. Otherwise
write $\ine(\alpha) = e_1 \ine(b_1) + \dots + e_r \ine(b_r)$,
where $b_1, \dots, b_r \in B$ and
$e_1, \dots, e_m$ are non-negative integers.
Then $\alpha$ and $b_1^{e_1} \dots b_r^{e_r}$
have the same leading exponent; thus for some $c \in k$,
\[ \alpha_1 = \alpha - c b_1^{e_1} \dots b_r^{e_r} \] has a lexicographically
smaller leading monomial than $\alpha$. If $\alpha_1 = 0$, we are done.
If not, we can replace $\alpha$ by $\alpha_1$ and apply the same procedure.
That is, after subtracting a monomial in elements of $B$ from $\alpha_1$,
we obtain $\alpha_2 \in R$ with a smaller initial exponent, etc.

In order to complete the proof of the proposition, it is enough to show
that the resulting sequence $\alpha = \alpha_0, \alpha_1, \alpha_2, \dots$
in $R$ will terminate, i.e.,  $\alpha_r = 0$ for some $r \geq 0$.
This is a very special case of~\cite[Proposition 6.5]{reichstein};
for the sake of completeness we give a direct proof below.

By our construction $\ine(\alpha_0) \succ \ine(\alpha_1) \succ
\ine(\alpha_2) \succ \dots$. Thus it suffices to prove that
for any given $a = (a_1, \dots, a_n) \in
\ine(R)$ there are only finitely many
$a' = (a_1', \dots, a_n') \in \ine(R)$ such that
$a \succeq a'$.
Indeed, if $a \succeq a'$ then $0 \leq a_1' \leq a_1$. Now condition
(iii) says that $a_n' \geq - a_1' \geq -a_1$, and condition (ii) says that
$-a_1 \leq a_n' \leq a_i' \leq a_1' \leq  a_1$. Thus
$a_i'$ may assume only finitely many values
for every $i = 1, \dots, n$. This completes the proof
of Proposition~\ref{prop.subduct}.
\end{proof}

For computational purposes, we are interested in those cases, where
the set $B$ in Proposition~\ref{prop.subduct} can be chosen to be finite,
i.e., $\ine(R)$ is a finitely generated semigroup. In such cases we shall
refer to $B$ as a SAGBI basis of $G$; cf.~\cite[Introduction]{reichstein}.
(Here SAGBI stands for ``subalgebra analog to Gr\"obner bases for ideals";
this term is due to Robbiano and Sweedler~\cite{rs}.)

Unfortunately, by~\cite[Theorem 1.6]{reichstein}
$k[X^{\pm 1}]_0^G$ has a SAGBI basis only for $n = 2$,
$3$ and $4$; see also \cite[Example 7.3]{reichstein}.
Moreover, the situation cannot be remedied by
replacing the lexicographic order with a different term order.
On the other hand, for $n \le 4$
the subduction algorithm is much faster
than Algorithm~\ref{alg.single}.

\subsection*{Explicit SAGBI bases}

Let $c_i$ be the $i$th elementary symmetric polynomial in
$x_i/x_j$, as in Lemma~\ref{lem2.1}.  (Here $i$ and $j$ are distinct
integers ranging from $1$ to $n$.)
Recall that $c_1 = q_1 -n = \tr(X) \tr(X^{-1})-n $.

\begin{lem} \label{lem.sagbi}
The following elements form a SAGBI basis of $k[X^{\pm 1}]_0^G$.

(a) $c_1$, if $n = 2$.

\smallskip
(b) $c_1$ and $c_2$, if $n = 3$.

\smallskip
(c) $c_1, c_2, c_3$ and $p$, if $n = 4$.
Here $p = s_2(X) s_2(X^{-1})$, where
\[ s_2(X) = x_1 x_2  + x_1 x_3 + \dots + x_3 x_4 \]
is the second symmetric polynomial in $X = (x_1, \dots, x_4)$.
\end{lem}

\begin{proof} Let $S = \ine(k[X^{\pm 1}]_0^G)$ be the subsemigroup
of $\bZ^n$ given by~\eqref{e.S}.

\smallskip
(a) If $n = 2$ then $S$ is clearly generated
by $(1, -1)  = \ine(c_1)$ as a semigroup.

\smallskip
(b) For $n = 3$, $S$ is generated, as a semigroup, by the elements
$\lambda_1 = (1, 0, -1) = \ine(c_1)$ and
$\lambda_2 = (2, -1, -1) = \ine(c_2)$.  Indeed,
every element of $\mu \in S$ is of the form $\mu = (a, -c, -b)$,
where $b \geq c \geq 0$ and $a = b + c$. Thus $\mu = c \lambda_1 +
(b-c) \lambda_2$ lies in the semigroup generated by $\lambda_1$ and
$\lambda_2$.

(c) We want to show that any
$\mu = (a, b, c, d) \in S$ can be written as a non-negative
integer linear combination of
\[ \begin{array}{l} \lambda_1 = (1, 0, 0, -1) = \ine(c_1) \, , \\
\lambda_2 = (2, 0, -1, -1) = \ine(c_2) \, , \\
\lambda_3 = (3, -1, -1, -1) = \ine(c_3) \, \; \; \text{and} \\
\lambda_4 = (1, 1, -1, -1) = \ine(p) \, . \end{array} \]
If $b \leq 0$ then the desired linear combination is given by
\[ \mu = (c-d)\lambda_1 + (b-c)\lambda_2 +
(-b ) \lambda_3 \, . \]
If $b > 0$ then, after replacing $\mu$ by $\mu - b\lambda_4$, we
can assume $b = 0$ and apply the above formula.
\end{proof}

\begin{cor} \label{cor.sagbi}
(a) If $n = 2$ then $k[X^{\pm 1}]_0^G = k[c_1] = k[q_1]$.

\smallskip
(b) If $n = 3$ then
$k[X^{\pm 1}]_0^G = k[c_1, \, c_2] = k[q_1, \, q_2]$.

\smallskip
(c) If $n = 4$, $k[X^{\pm 1}]_0^G = k[c_1, c_2, c_3, p] =
k[q_1, q_2, q_3, p]$.
\end{cor}

\begin{proof}
Immediate from Lemma~\ref{lem.sagbi}, Proposition~\ref{prop.subduct}
and Newton's formulas; cf. Lemma~\ref{lem2.1}(a).
\end{proof}

\begin{example} \label{ex.n=3} Let $n = 3$.
Corollary~\ref{cor.sagbi} tells us that
$c_1$ and $c_2$ form a SAGBI basis for $k[X^{\pm 1}]_0^G$.
For instance, applying the subduction algorithm, we obtain:
\[ E_2(X,X) = 2 (c_1^2 +  c_1 -  c_2) \, .\]
Similarly,
 \[ E_2(X,X^2) = 2 c_1^3 + 5 c_1^2 - 5 c_1 c_2 + 9 c_1 - 12 c_2 + 18 \, .\]
These identities will be used in Section~\ref{sect.explicit-phase}.
\end{example}

\begin{example} \label{ex.n=4-a} $n = 4$.
Using the subduction algorithm of Proposition~\ref{prop.subduct}
to express $c_2^2$ and $c_4$
in terms of the SAGBI basis
$c_1, c_2, c_3, p$, we obtain the following relations
in the ring $k[X^{\pm}]_0^G$:
\begin{eqnarray*}
  c_2^2 & = & 2 c_1^2 p - 16 c_1^2 - 8 c_1 c_2 - c_1 p^2 +
15 c_1 p - 48 c_1 + \\
   && 3 c_2 p - 12 c_2 + c_3 p - 2 p^2 + 18 p - 36  \\
c_4 & = & 3 c_1^2 + c_1 c_2 - 3 c_1 p + 17 c_1 + c_2 - 3 c_3 +
p^2 - 10 p + 21 \, . \\
\end{eqnarray*}
Eliminating $p^2$ and solving for $p$ gives
$$p = \dfrac{ 6 + 7c_1 + 7c_1^2
+ 3c_1^3 - 10c_2 - 5c_1c_2 + c_1^2c_2 - c_2^2
    - 6c_3 - 3c_1c_3 - 2c_4 - c_1c_4}{2+ c_1 + c_1^2 -3c_2 -c_3}.$$
This shows that
$p \in k(c_1, c_2, c_3, c_4)$ or, equivalently,
$p \in k(q_1, \dots, q_4)$. Thus by Corollary~\ref{cor.sagbi}
$k(X)_0^G = k(q_1, q_2, q_3, q_4)$. Recall that Theorem~\ref{thm1}
asserts only that
$k(X)_0^G = k(q_1, \dots , q_6)$; we have thus shown that the last
two of these generators are not needed.

To obtain explicit expressions for $c_5$ and $c_6$ as rational
functions in $c_1, \dots, c_4$, we use the subduction algorithm
once again:
\begin{eqnarray*}
  c_5 & = & c_1^3 - 6 c_1^2 - 5 c_1 c_2 + 7 c_1 p - 38 c_1 +
c_2 p - 7 c_2 + 6 c_3 - 2 p^2 + 20 p - 42 \\
  c_6 & = & -2 c_1^3 + c_1^2 p + 6 c_1 c_2 - 5 c_1 p + 27 c_1 - 2 c_2 p +
9 c_2 - 7 c_3 + \\ &&  2 p^2 - 18 p + 34 \, ,
\end{eqnarray*}
then substitute the above formula for $p$.
\end{example}

Note also that for $n = 3$,
Theorem~\ref{thm1} says that $k(X)_0^G = k(q_1, q_2, q_3)$, but
$q_3$ is not needed by Corollary~\ref{cor.sagbi}(b). We do not know
whether or not any of the generators listed in Theorem~\ref{thm1} can
be left out for $n \ge 5$.

\section{More sets of variables}
\label{sect.mult}

The purpose of this section is to generalize Theorem~\ref{thm1}
to the multi-array case. That is, instead of considering a single
array of independent variables $X= (x_1, \dots, x_n)$, we shall
consider $m$ arrays:
\[ \begin{matrix}
X_1 = (x_{11}, \dots, x_{1n}) \, , \\
X_2 = (x_{21}, \dots, x_{2n}) \, , \\
\hdotsfor{1} \\
X_m = (x_{m1}, \dots, x_{mn}) \, .
\end{matrix} \]
We shall view such $n$-tuples as diagonal $n \times n$-matrices
and operate with them as we did in Section~\ref{sect1}. For example,
\[  X_i X_j = (x_{i1} x_{j1}, \dots, x_{in} x_{jn}) \ , \quad
\tr(X_i) = x_{i1} + \dots + x_{in} \, , \]
the observables~\eqref{e.observables} can be written as
\begin{equation} \label{e.observables-trace}
q_{r_1, \dots, r_m}(X_1, \dots, X_m) =
\tr(X_1^{r_1} \dots X_m^{r_m}) \tr(X_1^{-r_1} \dots X_m^{-r_m})
\end{equation}
and the functions $E_m$ that arise in the phase transition
problem~\eqref{e.E_m} as
\begin{eqnarray} \label{e.E_m-trace} E_m(X_1, \dots, X_m) =
\tr(X_1) \dots \tr(X_m) \tr(X_1^{-1} \dots X_m^{-1}) + \\ \nonumber
\tr(X_1^{-1}) \dots \tr(X_m^{-1}) \tr(X_1 \dots X_m) \, .
\end{eqnarray}

Let $k(X_1, \dots, X_m)_0$
be the subfield of $k(x_{11}, x_{12}, \dots, x_{mn})$
whose elements are rational functions in $x_{ij}$ homogeneous
of degree 0 in each $n$-tuple of variables $X_i$. In other words,
$k(X_1, \dots, X_m)_0$
is the function field of the variety $(\bP^{n-1})^m$.
In Section~\ref{sect1} we studied the action of the group
$G = \Sym_n \times \T$ on $\bP^{n-1}$;
this action extends to an action of $G \times \dots \times G = G^m$
on $(\bP^{n-1})^m$. We shall be interested in the invariants
for the action of the diagonal subgroup of $G^m$
which we shall also denote by $G$.  In concrete terms,
the symmetric group $\Sym_n$ acts on the function field
$k(X_1, \dots, X_m)_0$ of $(\bP^{n-1})^m$ by simultaneously
permuting the variables $x_{i1}, \dots, x_{in}$ for each $i$.
The 2-element group $\T = \{ 1, \tau \}$ acts on
$k(X_1, \dots, X_m)_0$ by
\begin{equation} \label{e.tau-mult}
\text{$\tau \colon x_{ij} \lra \frac{1}{x_{ij}}$ for every $i = 1, \dots, m$
and $j = 1, \dots, n$.}
\end{equation}
These two actions commute and thus induce an action of
$G = \Sym_n \times \T$ on $k(X_1, \dots, X_m)_0$.

\begin{thm} \label{thm2} The field $k(X_1, \dots, X_m)_0^G$
is generated (over $k$) by the following elements:
$$ q_{1, \dots, 1} := \tr(X_1 \dots X_m) \tr \, (X_1 \dots X_m)^{-1}$$
and
\[ \begin{array}{c}
q_{r, 0, \dots, 0} := \tr(X_1^r) \tr(X_1^{-r}) \, , \\
q_{0, r, \dots, 0} := \tr(X_2^r) \tr(X_2^{-r})\, , \\
 \vdots \\
q_{0, \dots, 0, r} := \tr(X_m^r) \tr(X_m^{-r})\, ,
\end{array} \]
where $r$ ranges from $1$ to $(n(n-1))/2$.
\end{thm}

We will give a proof that generalizes the proof in Section~\ref{sect1},
and then discuss the prospects for a more constructive proof.
Note that our earlier results for $m =1$ will be
used in the proof.

We begin by disposing of the case $n = 2$. This case is anomalous in that
$G$ does not act effectively on
$k(X_1, \dots, X_m)_0$; the kernel of this action is the 2-element
subgroup of $G = \Sym_2 \times \T$ generated by $(\sigma, \tau)$,
where $\sigma$ is the non-trivial element of $\Sym_2$ and
$\tau$ is the non-trivial element of $\T$. Thus for $n = 2$
\[ k(X_1, \dots, X_m)_0^G = k(X_1, \dots, X_m)_0^{\tau} \, ; \]
here $\tau$ acts via the involution~\eqref{e.tau-mult}. Setting
\[ t_i = \frac{x_{i1}}{x_{i2}} \]
we see that $k(X_1, \dots, X_m)_0 = k(t_1, \dots, t_m)$ and
$\tau$ acts on this field by taking $t_i$ to $\frac{1}{t_i}$ for
each $i = 1, \dots, m$. It is now a simple exercise in Galois
theory to show that in this case
\[ k(t_1, \dots, t_m)^{\tau} = k(t_1 + \frac{1}{t_1}, \dots,
t_m + \frac{1}{t_m}, t_1 \dots t_m + \frac{1}{t_1 \dots t_m}) \, . \]
Since $t_i + \frac{1}{t_i} = \tr(X_i)\tr(X_i^{-1}) - 2$ and
\[ t_1 \dots t_m + \frac{1}{t_1 \dots t_m} =
\tr(X_1 \dots X_m)\tr \, (X_1 \dots X_m)^{-1} -2 \, , \]
this proves Theorem~\ref{thm2} for $n = 2$.

{}From now on we will assume that $n \geq 3$.
Let $K$ be the field generated over~$k$ by the observables
listed in the theorem.  For $h = 1, \dots, n$ let
\begin{equation} \label{e.f_i}
f_h(t) = \prod_{i \neq j} (t-\frac{x_{hi}}{x_{hj}}) \, .
\end{equation}
By Lemma~\ref{lem2.1} the coefficients of each $f_h(t)$ lie in $K$.
The field $k(X_1, \dots, X_m)_0$ is clearly the splitting field
of the product \[ f(t) := f_1(t) \dots f_m(t) \] over $K$.
The Galois group of~$f$ contains~$G$ and is naturally embedded in the
product $G^m$ of the Galois groups of the $f_i$.

\begin{lem} \label{lem3.6} Assume $n \geq 3$.
An element $g = (g_1, \dots, g_m)$ of $G^m$
fixes \[ q_{1, \dots, 1} = \tr(X_1 \dots X_m)
\tr \, (X_1 \dots X_m)^{-1} \in k(X_1, \dots, X_m)_0^G \]
if and only if $g_1 = \dots = g_m$.
\end{lem}

\begin{proof} One direction is obvious: $(g_1, \dots, g_1)$
fixes $q_{1, \dots, 1}$ for
every $g_1 \in G$.  For the purpose of proving the converse,
we may replace $g$ by $(g_1^{-1}, \dots, g_1^{-1}) g$ and thus assume
$g_1 = id$. In other words, we want to prove that if
$g = (id, g_2, \dots, g_n)$ fixes $q_{1, \dots, 1}$
then $g_2 = \dots g_m = id$ in $G$.

For $i = 2, \dots, m$, let $g_i = (\sigma_i, \epsilon_i) \in G
= \Sym_n \times \T$.
Here $\T = \{ 1, \tau \}$ is written multiplicatively, i.e.,
$\tau$ is written as $-1$ and each $\epsilon_i = \pm 1$.
In particular,
\[ (1, g_2, \dots, g_m) \cdot (\frac{x_{hi}}{x_{hj}}) =
(\frac{x_{h\sigma_h(i)}}{x_{h \sigma_h(j)}})^{\epsilon_h} \, . \]
Writing out $q_{1, \dots, 1}$ explicitly in terms of the $x_{ij}$
(cf.~\eqref{e.observables}), we obtain
\[ q_{1, \dots, 1} -n =
\sum_{i \neq j} \frac{x_{1i}}{x_{1j}} \frac{x_{2i}}{x_{2j}}
 \dots \frac{x_{mi}}{x_{mj}} \in K \, ,  \]
and
\[ \label{e.g-tr(xy)}
q_{1, \dots, 1} -n = g(q_{1, \dots, 1} - n) =
\sum_{i \neq j} \frac{x_{1i}}{x_{1j}}
(\frac{x_{2\sigma_2(i)}}{x_{2\sigma_2(j)}})^{\epsilon_2}
\dots
(\frac{x_{m \sigma_m(i)}}{x_{m \sigma_m(j)}})^{\epsilon_m} \, .
\]
Comparing the terms of the last two equations, we see that for each
$h = 2, \dots, m$, either (i) $\epsilon_h = 1$ and $\sigma_h(i) = i$
for every $i = 1, \dots, n$, i.e., $\sigma_h = id$, or
(ii) $\epsilon_h = -1$ and
$\sigma_h(i) = j$ for every pair of distinct integers $i, j$ between
$1$ and $n$.  In case (ii), $\sigma(i)$ assumes every value between
$1$ and $n$ other than $i$, which is impossible for $n \geq 3$.
We conclude that $g_2 = \dots = g_n$, as claimed.
\end{proof}

This shows that $G$ is the Galois group of
$k(X_1, \cdots , X_m)_0$ over~$K$. Examining the tower
\[ \begin{array}{c}   k(X_1, \dots, X_m)_0 \\
                          |   \\
                     k(X_1, \dots, X_m)_0^G \\
                         |   \\
                         K
\end{array} \]
we conclude that $k(X_1, \dots, X_m)_0^G = K$,
thus completing the proof of the theorem.

We now remark that Algorithm~\ref{alg.single} (for $m = 1$)
can be extended to the multi-array case by means of a suitable
universal denominator.  The generalization is not
entirely straightforward  because of the ``new" generator
\[ q_{1, \dots, 1} = \tr(X_1 \dots X_m)\tr \, (X_1 \dots X_m)^{-1} \]
that connects the different $X_i$.

Let $\Ri$ be the $k$-linear
combinations of Laurent monomials $x_{11}^{a_{11}} \dots x_{mn}^{a_{mn}}$
such that $a_{i1} + \dots + a_{in} = 0$ for every $i = 1, \dots, m$.
Note that $\Ri$ is generated,
as a $k$-algebra, by elements of the form $\frac{x_{ij}}{x_{il}}$;
in particular,
the field $k(X_1, \dots, X_m)_0$ is the field of fractions of $\Ri$.
Note also that $\Ri$ is a $G$-invariant subring of $k(X_1, \dots, X_m)_0$.

Let $z_{hij}$ be a set of $mN$ algebraically independent
variables over $k$, where $N = n(n-1)$, $i$ and $j$ are
distinct integers between
$1$ and $n$, and $h$ ranges from $1$ to $m$.
Let $t$ be another independent variable
and let
\[ \phi \colon k[z_{hij}, t] \lra \Ri \]
be the surjective $k$-algebra homomorphism
given by
\[ \text{$\phi(z_{hij}) = \frac{x_{hi}}{x_{hj}}$
and $\phi(t) = q_{1, \dots, 1} - n$.} \]
We note that $\Sym_N^m$ acts on $k[z_{hij},t]$ and
where the $h$-th component permutes the $z_{hij}$ and fixes~$t$.

Our universal denominator is a polynomial $E \in k[z_{hij},t]$
that has the property that for $\sigma_1, \dots, \sigma_m \in \Sym_{N}$
\begin{equation} \label{e.E}
\phi((\sigma_1, \dots, \sigma_m)E) = \begin{cases} \phi(E) &
\text{if $\sigma_1 = \dots = \sigma_m \in G$,} \\
0 & \text{in all other cases.}
\end{cases}
\end{equation}
Here, as before, we view $G$ as a subgroup of $\Sym_N$ for $n \geq 3$.

The polynomial $E$ is defined by
$$E(z_{hij},t) := D(z_{1ij}) \cdot \cdots \cdot D(z_{mij}) \,
E_1(z_{hij},t) \, , $$
where $D$ is the polynomial defined in Section~\ref{sect.alg-single} and
$$E_1(z_{hij},t) := \prod
 (t - \sum_{i \ne j} z_{1ij} z_{2 g_2(ij)}
\dots z_{m g_m(ij)} ) \, ;$$
the product is taken over all $g_2, \dots, g_m \in G$ such that
at least one $g_i \neq 1$, and
the sum is taken over all pairs of distinct integers
$i$ and $j$ between $1$ and $n$.

As in the $m = 1$ case one can verify that $E$
satisfies~\eqref{e.E}.  The algorithm for expressing invariants
in terms of observables is now similar to the $m = 1$ case,
and we leave the details to the reader.

\section{Reduction to one set of variables}
\label{sect.explicit-phase}

We now return to the problem of writing a $G$-invariant multihomogeneous
rational function $f(x_{11}, x_{12}, \dots, x_{mn})$
of total degree 0 in each group of variables
$X_1 = (x_{11}, \dots, x_{1n})$, $\dots$, $X_m = (x_{m1}, \dots, x_{mn})$,
as a rational function in the ``observables" $q_{r_1, \dots, r_m}$.

In principle, the algorithm sketched in the preceding section is a
solution to this problem.  As one might expect, it is too slow
to be of practical significance. On the other hand,
the approach we took in Section~\ref{sect.sagbi}
cannot be extended to $m \geq 2$,
even for small $n$, because suitable (finite) SAGBI bases
do not exist; see~\cite[Theorem 1.6 and Example 7.3]{reichstein}.
This motivated our search for an algorithm that would reduce computations
in the multi-array case ($m \geq 2$) to computations in
the single-array case ($m = 1$). In this section we discuss such
an algorithm and
use it to generate explicit formulas for $E_2(X, Y)$ for $n = 3$ and $4$.

\subsection*{Another generating set of observables}
We begin by proving another variant of Theorem~\ref{thm0.1}(a).

\begin{thm} \label{thm3} The field $k(X_1, \dots, X_m)_0^G$
is generated (over $k$) by the elements
$$ q_r := q_r(X_1) = q_{r, 0, \dots, 0} := \tr(X_1^r) \tr(X_1^{-r}) $$
and
\[ \begin{array}{c}
q_{s}^{(2)} := q_{s, 1, 0, \dots, 0} := \tr(X_1^s X_2) \tr(X_1^{-s}X_2^{-1}) \, , \\
 \vdots \\
q_{s}^{(m)} := q_{s, 0 \dots, 0, 1} := \tr(X_1^s X_m) \tr(X_1^{-s} X_m^{-1})
\, .  \end{array} \]
where $r = 1, \dots, \frac{n(n-1)}{2}$ and $s = 0, \dots, n(n-1) - 1$.
\end{thm}

Informally speaking, the element $q_{1, \dots, 1}$ of Theorem~\ref{thm2}
ties $X_1, \dots, X_m$ together, where as the elements $q_{r}^{(i)}$
of Theorem~\ref{thm3} relate $X_i$ to $X_1$ for each $i = 2, \dots, m$.

\begin{proof}
Let $K = k(X_1, \dots, X_m)_0 = k(x_{ij}/x_{il})$,
and let $K^N$ be a vector space of dimension $N = n(n-1)$ over $K$.
We shall write elements of $K^N$ as $(z_{ij})$, where $(i, j) \in \Lambda$,
i.e., $i$ and $j$ are distinct integers between $1$ and $n$,
as in~\eqref{e.Lambda}. The natural action of $G = \Sym_n \times T$
on $\Lambda$ induces a permutation action on $K^N$.

Given an $n$-tuple $A = (a_1, \dots, a_n)$, we will denote the
$N$-tuple of ratios $\frac{a_i}{a_j}$ by $\rho(A)$. Then
$\rho \colon K^n \lra K^N$ is a $G$-equivariant map.
Finally, for $h = 1, \dots, m$ let $Z_h$ be the $N$-tuple
of variables $(z_{h ij})$, where $(i, j) \in \Lambda$.

We want to show that any $f(X_1, \dots, X_m) \in k(X_1, \dots, X_m)^G$
can be written as a rational function in the observables listed in the
statement of the theorem. We begin with several reductions.
First of all, we may assume without loss of generality
that $f(X_1, \dots, X_m) \in \Ri^G$, by writing an invariant
rational function as a quotient of invariant polynomials.
Secondly, $f(X_1, \dots, X_m)$ can be lifted
to a $G$-invariant polynomial $F(Z_1, \dots, Z_m)$ in $z_{hij}$
such that \[ f(X_1, \dots, X_m) = F(\rho(X_1),
\dots, \rho(X_m)) \, .   \]
Thirdly, we may assume without loss of generality
that $F(Z_1, \dots, Z_m)$ is a homogeneous polynomial in the
arrays of variables $Z_1, \dots, Z_m$ of multi-degree
$(d_1, \dots, d_m)$. Indeed, any $G$-invariant $F$
can be written as a sum of $G$-invariant multihomogeneous components,
say, $F = F_1 + \dots + F_r$, and we may replace $f$ by
$f_i = F_i(\rho(X_1), \dots, \rho(X_m))$.

Multilinearizing $F$, we obtain a $G$-invariant multilinear
polynomial $M$ in $d = d_1 + \dots + d_m$
$N$-variable arrays such that
\begin{equation} \label{e.multilinear}
f(X_1, \dots, X_m) = M(
\underbrace{\rho(X_1), \dots, \rho(X_1)}_{\text{$d_1$ times}}, \dots,
\underbrace{\rho(X_m), \dots, \rho(X_m)}_{\text{$d_m$ times}}) \, .
\end{equation}

Next we observe that by the Vandermonde argument
\begin{equation} \label{e.basis0}
\rho(I), \rho(X_1), \dots, \rho(X_1^{N-1})
\end{equation}
form a $K$-basis of $K^N$; here $I$ stands for the identity $n$-tuple
$(1, \dots, 1)$. In particular, for every $2 \le i \le m$, we can write
\begin{equation} \label{e0}
\rho (X_i) = \lambda_{i0} \rho(I) + \lambda_{i1} \rho(X_1) + \dots +
\lambda_{i, N-1} \rho(X_1^{N-1}) \, .
\end{equation}
for some $\lambda_{i0}, \dots, \lambda_{i \, N-1} \in K^N$.
Substituting this into~\eqref{e.multilinear} and expanding, we see
that $f(X_1, \dots, X_m)$ can be written as a sum of terms of the form
\[ \text{(monomial in $\lambda_{ij}$)
 $M(\rho(X_1^{i_1}), \dots, \rho(X_1^{i_d}))$.} \]
(In fact, $i_1 = \dots = i_{d_1} = 1$ in each term, but we shall not
use this in the sequel.)
Since $M$ is $G$-invariant, each
\[ M(\rho(X_1^{i_1}), \dots, \rho(X_1^{i_d})) \]
is an element of $k[X_1]_0^G$, and thus, by Theorem~\ref{thm1}, it
can be written as a rational function in the observables $q_r = q_r(X_1)$.

Thus it remains to show that each $\lambda_{ij}$ lies in the field $L$
generated by the elements listed in the statement of the theorem. Note that
by Theorem~\ref{thm1} the observables $q_r$ lie in $L$ for every $r$
(and not just for $r = 1, \dots N/2$).
Taking the dot product of both sides of~\eqref{e0}
with $\rho(I), \rho(X_1), \dots, \rho(X_1^{N-1})$, and remembering that
\[ \text{$ \rho(X_1^i) \cdot \rho(X_1^j) =
q_{i+j} - n$ and $\rho(X_i) \cdot \rho(X_1^j) = q_{j}^{(i)} - n$,} \]
we obtain the following linear system:
\begin{equation} \label{e.linear-system}
\text{$\sum_{i=1}^N (q_{i+j} - n) \lambda_i =
 q_{j}^{(i)} - n$, where $i = 0, 1, \dots, N-1$.}
\end{equation}
Since $\rho(I), \rho(X_1), \dots, \rho(X_1^{N-1})$
form a basis of $K^N$, the matrix
\begin{equation} \label{e.matrixQ} Q = \begin{pmatrix}
q_0-n & q_{1} - n & \dots & q_{N-1} - n \\
q_{1} -n & q_{2} - n & \dots & q_{N-1} - n \\
\dots \\
q_{N-1} -n & q_{N} - n & \dots & q_{2N-2} - n \end{pmatrix} \, .
\end{equation}
of this system is non-singular.
Solving the linear system~\eqref{e.linear-system},
we conclude that each $\lambda_j$ lies in $L$, as claimed.
\end{proof}

\subsection*{Multilinear invariants}

Our proof of Theorem~\ref{thm3} reduces the computation of a given
element
$f$ of $k(X_1, \dots, X_m)_0^G$ to computing finitely many
elements in $k[X_1]_0^G$. Note that this reduction is constructive.
The resulting algorithm is cumbersome in general; however, it simplifies
considerably in the case where the polynomial $F(Z_1, \dots, Z_m)$
is itself multilinear, i.e., $d_1 = \dots = d_m = 1$
and $M = F$. This is exactly what happens
in that case of greatest interest to us, namely,
$f = E_m$ (cf.~\eqref{e.E_m}); here
\[ F(Z_1, \dots, Z_m) = \sum_{j_1, \dots, j_m, j = 1}^n
z_{1 j_1 j} \dots z_{m j_m j} +
z_{1 j j_1} \dots z_{m j j_m} \, . \]

\begin{prop} \label{prop.multilinear}

\[E_m(X_1, \dots, X_m) = \sum_{i_2, \dots, i_m = 0}^{N-1}
\lambda_{2i_1} \dots \lambda_{m i_m} E_m(X_1, X_1^{i_2}, \dots, X_1^{i_m})
\, , \]
where
\[ \begin{pmatrix}
\lambda_{i0} \\ \lambda_{i1} \\ \dots \\ \lambda_{i \, N-1}
\end{pmatrix} = Q^{-1}
\begin{pmatrix}
q_0^{(i)} - n \\ q_1^{(i)} - n \\ \dots \\ q_{N-1}^{(i)} - n
\end{pmatrix} \, . \]
and $Q$ is the $N \times N$-matrix~\eqref{e.matrixQ}.
\end{prop}

\begin{proof} The first formula is obtained by substituting~\eqref{e0}
into
\[ E_m(X_1, \dots, X_m) = F(\rho(X_1), \dots, \rho(X_m)) \]
and expanding the right hand side. (Note that the specific form of $F$
is not used here; we only use the fact that $F$ is multilinear.)
The second formula comes from solving
the linear system~\eqref{e.linear-system}.
\end{proof}

\begin{remark} Proposition~\ref{prop.multilinear} remains valid
if $E_m$ is replaced by any $f$ in $ \Ri^G$
such that
$f(X_1, \dots, X_m) = F(\rho(X_1), \dots, \rho(X_m))$
for some $G$-invariant
multilinear polynomial $F(Z_1, \dots, Z_m)$ in $m$ $N$-variable
arrays $Z_1, \dots, Z_m$.
\end{remark}

\subsection*{Explicit determination of the triplet phase invariant}

For $m = 2$ the formula of Proposition~\ref{prop.multilinear} can be
rewritten in the matrix form:
\begin{equation} \label{e.explicit-formula1}
E_2(X, Y) =  (e_0, \dots , e_{N-1}) \, Q^{-1}
\begin{pmatrix}
q_{0, 1} - n \\ q_{1, 1} - n \\ \dots \\ q_{N-1, 1} - n
\end{pmatrix} \, , \end{equation}
where $e_i = e_i(X) = E_2(X, X^i)$.
Here, as usual, we write $X$ for $X_1$ and $Y$ for $X_2$; cf.~\eqref{e.E_2}.
We have thus reduced the computation of $E_2(X, Y)$ to the computation of
$e_i = E(X, X^i)$ for $i = 0, \dots, N-1$. Note that
$e_0 = 2nq_1$, so there are only
$N-1$ elements $e_1, \dots, e_{N-1}$ we need to compute.
We can further reduce this number by using the basis
\[ \rho(X^{-N/2}), \rho(X^{-N/2 + 1}), \dots, \rho(X^{N/2 -1}) \]
of $K^n$ instead of~\eqref{e.basis0}. This has the effect of shifting the
subscripts in~\eqref{e.explicit-formula1} as follows:

\begin{prop} \label{prop.explicit-formula2}
\[ E_2(X, Y) = (e_{-N/2}, e_{1-N/2}, \dots, e_{N/2 -1}) \cdot R^{-1} \cdot
\begin{pmatrix}
q_{-N/2, 1} - n \\ q_{-N/2 + 1, 1} - n \\ \dots \\ q_{N/2 -1, 1} - n
\end{pmatrix} \, , \]
where $q_r = q_{r, 0}$, $e_i = E_2(X, X^i)$, and $R$ is
the $N \times N$-matrix
\[ R = \begin{pmatrix}
q_{-N}-n & q_{-N + 1} - n & \dots & q_{-1} - n \\
q_{-N + 1} -n & q_{-N + 2} - n & \dots & q_0 - n \\
\dots \\
q_{-1} -n & q_{0} - n & \dots & q_{N-2} - n \end{pmatrix} \, . \]
\qed
\end{prop}
Keeping in mind the identities $e_i = e_{-1-i}$
and $e_0 = 2nq_1$, we see that
the formula of Proposition~\ref{prop.explicit-formula2} reduces the
computation of $E_2(X, Y)$ to the computation of $e_i = E_2(X, X^i)$ for
$i = 1, \dots, (N/2)-1$. We also note that since $q_{-r} = q_r$,
the matrix $R$ involves $q_r$
only for $r = 1, \dots, N$, as opposed to $r = 1, \dots, 2N-2$
for the matrix~\eqref{e.matrixQ}.

The calculation of the $e_i$ involves only one set of variables.
For $n = 3$ we have the following explicit results.

\begin{example} \label{ex.formula3}
Let $n = 3$. Then $N = 6$, $q_0 = 9$,
and Proposition~\ref{prop.explicit-formula2} gives
\[ E_2(X, Y) = (e_2, e_1, 6q_1, 6q_1, e_1, e_2) \, R^{-1} \,
\begin{pmatrix}
q_{-3, 1} - 3 \\ q_{-2, 1} - 3 \\ \dots \\ q_{2, 1} - 3
\end{pmatrix} \, , \]
where $R$ is the $6 \times 6$-matrix
\[ R = \begin{pmatrix}  q_6-3 & q_5-3 & q_4-3 & q_3-3 & q_2-3 & q_1-3 \\
q_5-3 & q_4-3 & q_3-3 & q_2-3 & q_1-3 & 6 \\
q_4-3 & q_3-3 & q_2-3 & q_1-3 & 6 & q_1-3 \\
q_3-3 & q_2-3 & q_1-3 & 6 & q_1-3 & q_2-3 \\
q_2-3 & q_1-3 & 6 & q_1-3 & q_2-3 & q_3-3 \\
q_1-3 & 6 & q_1-3 & q_2-3 & q_3-3 & q_4-3 \end{pmatrix} \, .\]
Recall that $e_1$ and $e_2$ were computed in Example~\ref{ex.n=3}
by using SAGBI basis techniques:
\begin{eqnarray*}
  e_1 & = & 2 c_1^2 + 2 c_1 - 2 c_2 \\
  e_2 & = & 2 c_1^3 + 5 c_1^2 - 5 c_1 c_2 + 9 c_1 - 12 c_2 + 18 \, .\\
\end{eqnarray*}
These can be easily expressed in terms of $q_1$ and $q_2$ by using
Newton's identities, which here amount to
$$ c_1 = q_1 - 3 \qquad \text{and} \qquad  c_2 = ((q_1 - 3)^2 - (q_2 - 3))/2.$$
The explicit formula for $E_2$ that results is
quite different from,
and considerably more elaborate than, the formula given in
Example~\ref{first.threeformula}.
\end{example}

\begin{example} \label{ex.formula4}
Let $n = 4$. Then $N = 12$, $q_0 = 16$,
and the formula of Proposition~\ref{prop.explicit-formula2} reduces to
\[ E_2(X, Y) = (e_5, e_4, e_3, e_2, e_1, 8q_1, 8q_1, e_1, e_2, e_3, e_4, e_5)
\, \cdot R^{-1} \, \cdot
\begin{pmatrix}
q_{-6, 1} - 4 \\ q_{-2, 1} - 4 \\ \dots \\ q_{5, 1} - 4
\end{pmatrix} \, , \]
where $R$ is the $12 \times 12$-matrix
\[ R = \begin{pmatrix}  q_{12}-4 & q_{11}-4 & \dots & q_3-4 & q_2-4 & q_1-4 \\
q_{11}-4 & q_{10}-4 & \dots & q_2-4 & q_1-4 & 12 \\
q_{10}-4 & q_9-4 & \dots & q_1-4 & 12 & q_1-4 \\
\cdots & \cdots & \cdots & \cdots & \cdots & \cdots \\
q_1-4 & 12 & \dots & q_{9}-4 & q_{10}-4 & q_{11}-4 \end{pmatrix} \, .\]
The elements $e_1, \dots, e_5$ can again be explicitly
determined using the SAGBI basis subduction algorithm; the
result is:

\begin{eqnarray*}
  e_1 & = & 2 c_1^2 + 8 c_1 - 2 c_2 - 4 p + 20\\
  e_2 & = & 2 c_1^3 + 4 c_1^2 - 5 c_1 c_2 - c_1 p + 7 c_1 - 5 c_2 + 3 c_3 + 2 p + 2\\
  e_3 & = & 2 c_1^4 + 4 c_1^3 - 7 c_1^2 c_2 + 3 c_1^2 p - 41 c_1^2 - 31 c_1 c_2  \\
       && + 7 c_1 c_3 - 2 c_1 p^2 + 40 c_1 p - 154 c_1 + 8 c_2 p - 28 c_2  \\
       && + 2 c_3 p + 24 c_3 - 8 p^2 + 72 p - 136 \\
  e_4 & = & 2 c_1^5 + 4 c_1^4 - 9 c_1^3 c_2 + 13 c_1^3 p - 131 c_1^3 - 80 c_1^2 c_2
      + 9 c_1^2 c_3 \\
      &&  - 7 c_1^2 p^2 + 144 c_1^2 p - 651 c_1^2 + 24 c_1 c_2 p - 186 c_1 c_2
         + 7 c_1 c_3 p  \\
      && + 43 c_1 c_3 - 29 c_1 p^2 + 379 c_1 p - 1136 c_1 - 5 c_2 c_3 + 21 c_2 p \\
      &&- 120 c_2 + 3 c_3 p + 78 c_3 - 34 p^2 + 334 p - 688\\
  e_5 & = & 2 c_1^6 + 4 c_1^5 - 11 c_1^4 c_2 + 27 c_1^4 p - 243 c_1^4 - 142 c_1^3 c_2 \\
       && + 11 c_1^3 c_3 - 14 c_1^3 p^2 + 289 c_1^3 p - 1317 c_1^3 + 46 c_1^2 c_2 p \\
       && - 429 c_1^2 c_2 + 14 c_1^2 c_3 p + 53 c_1^2 c_3 - 66 c_1^2 p^2 + 892 c_1^2 p\\
       && - 2830 c_1^2 - 17 c_1 c_2 c_3 + 81 c_1 c_2 p - 471 c_1 c_2 + 19 c_1 c_3 p \\
       && + 135 c_1 c_3 + 2 c_1 p^3 - 137 c_1 p^2 + 1286 c_1 p - 3039 c_1 - 2 c_2^3 \\
       && - 35 c_2 c_3     + 27 c_2 p - 167 c_2 + 3 c_3^2 - 2 c_3 p^2 - c_3 p \\
       && + 171 c_3 + 8 p^3 - 150 p^2 + 856 p - 1402
\end{eqnarray*}
Then we eliminate $p$, by using the formula in Example~\ref{ex.n=4-a}.
to express $e_1, \dots, e_5$ as rational functions in
$c_1$, $c_2$ and $c_3$. Finally, we use Newton's identities to
rewrite each $e_i$ as a rational function of $q_1$, $q_2$ and $q_3$.
\end{example}

\begin{remark}
The SAGBI basis algorithms were implemented explicitly in magma~\cite{magma},
and are very efficient, though of course they are limited to $n \le 4$.
\end{remark}

\section{Regular invariants}
\label{sect.regular}

Theorems~\ref{thm2} and~\ref{thm3} tell us that if we allow
$r_1, \dots, r_m$ to range over the integers then
the observables $q_{r_1, \dots, r_m}(X_1, \dots, X_m)$
generate $k(X_1, \dots, X_m)_0^G$ as a field
extension of $k$.
H.~Hauptman asked us whether the
functions $E_m$ in~\eqref{e.E_m-trace}
that arise in the phase determination problem~\eqref{e.E_m} can be
written as {\em polynomials}, rather than rational functions
in the $q_{r_1, \dots, r_m}$.
More generally, we can ask whether the the observables $q_{r_1, \dots, r_m}$
generate the $k$-algebra $\Ri^G$.
In this section, we will show that the answer is generally ``no.''

\begin{prop} \label{prop.n>=5} Assume $n \geq 4$.

\smallskip
(a) The function
\[ f(X) = E_2(X, X) = \tr^2(X) \tr(X^{-2}) + \tr^2(X^{-1}) \tr(X^2) \]
cannot be written as a polynomial in the observables $q_r$, as $r$
ranges over the integers.

\smallskip
(b) Suppose $m \geq 2$. Then function $E_m(X_1, \dots, X_m)$ given
by~\eqref{e.E_m-trace} cannot  be written as a polynomial
in the observables $q_{r_1, \dots, r_m}$,
as $r_1, \dots, r_m$ range over the integers.

\smallskip
(c) For any $m \geq 1$
the $k$-algebra $\Ri^G$ is not generated by the observables
$q_{i_1, \dots, i_m}$, as $i_1, \dots, i_m$ range over the integers.
\end{prop}

By Corollary~\ref{cor.sagbi}(a) and (b), $f$ can be written as a
polynomial in the observables if $n = 2$ or $3$; i.e., part (a) is
not true for $n = 2$ or~$3$.  Also, part (b) is not true for
$m = 1$ (indeed, $E_1$ is itself an observable).

Curiously, as we saw in the Section~\ref{sect.mult},
the observables $q_{i_1, \dots, i_m}$ come close to generating
the $k$-algebra $\Ri^G$ in the following sense:
every $\alpha \in \Ri^G$ can be written in the form
\[ \alpha = \frac{\beta}{\phi(E)} \] for some
$\beta \in k[q_{i_1, \dots, i_m}]$. Here $\phi(E)$ is a fixed
nonzero element of $k[q_{i_1, \dots, i_m}]$, independent of~$\alpha$.

\begin{proof} Part (b) is easily deduced from part (a) by specializing
$X_1$ and $X_2$ to a single array of indeterminates
$X = (x_1, \dots, x_n)$, and $X_3, \dots, X_m$ to the ``identity array"
$I = (1, \dots, 1)$. Part (c) is a consequence of parts (a)
(for $m = 1$) and (b) (for $m \geq 2$), since $f \in k[X]_0^G$ and
$E_m \in \Ri^G$ for every $m \ge 2$.

Thus we only need to prove (a).  Let
\[ \rho \colon \bA^n - \{ x_1 \dots x_n = 0 \} \lra \bA^{N} \]
be the map given by
\begin{equation} \label{e.mapf}
\rho(x_1, \dots , x_n) = (\frac{x_1}{x_2}, \frac{x_1}{x_3}, \dots,
\frac{x_{n-1}}{x_n}) \, . \end{equation}
To prove part (a), we claim that it suffices to exhibit $n$-tuples $x$ and $y$
with nonzero coordinates such that
\[ \text{(i) $f(x) \ne f(y)$  $\quad$ and $\quad$
(ii) $\rho(x) \sim \rho(y)$,} \]
where two $N$-tuples $u$ and $v$
are equivalent, written $u \sim v$, if one is a permutation of the other.

Indeed, assume that (i) and (ii) are true.
Since the observables $q_r$ are symmetric functions
in $\{ x_i/x_j \}$, where $i, j = 1, \dots, n$,
$i \neq j$, the fact that $\rho(x) \sim \rho(y)$ implies that
$q_r(x) = q_r(y)$ for every integer $r$.  Thus $q(x) = q(y)$
for every $q \in k[q_1, q_2, \dots]$; given that $f(x) \ne f(y)$,
we immediately deduce that
$f \not \in k[q_1, q_2, \dots]$ as desired.  Note that it is
sufficient if the coordinates of $x$ and $y$ lie in the
algebraic closure of the field~$k$.

First consider the case $n > 4$.  Let $z$ be a primitive $n$-th
root of unity.  (Recall that we are assuming that $\Char(k) = 0$,
so that $z$ exists in the algebraic closure of~$k$.) We claim that
the points
\begin{eqnarray*}
x & = & (1,1, z^3, z^3, z^4, \cdots , z^{n-1})\\
y & = & (z, z, z^2, z^2, z^4, \cdots , z^{n-1})
\end{eqnarray*}
satisfy conditions (i) and (ii).
Here $x$ is obtained from the point $(1, z, z^2,  \cdots, z^{n-1})$
by replacing $z$ by~$1$, and $z^2$ by $z^3$, and $y$ is obtained from
$(1, z , z^2,  \cdots, z^{n-1})$ by a similar alteration of two coordinates.
Note that
\begin{equation} \label{e.rt-of-unity}
\tr(x^i) + \tr(y^i) = 2(1^i + z^i + z^{2i} + \dots + z^{(n-1)i}) = 0
\end{equation}
for any $i \not \equiv 0 \pmod{n}$. In particular, $\tr(x^i) = - \tr(y^i)$ for
$i = \pm 1, \pm 2$ and consequently, $f(x) = - f(y)$. Thus, in order to check
(i), we only need to verify that $f(x) \ne 0$.
This is easily done; substituting
\[ \tr(x^i) = 1-z^i-z^{2i}+z^{3i} = (1-z^i)(1-z^{2i}) \]
into $f(x) = \tr^2(x) \tr(x^{-2}) + \tr^2(x^{-1}) \tr(x^2)$, we see that
\[
f(x) = 2 z^{-6} (1-z)^2(1-z^2)^3(1-z^4) \ne 0
\]
for any $n > 4$.

Now we have to show that $\rho(x) \sim \rho(y)$.  To do this
let $G_x(t)$ denote the formal generating function
$$G_x(t) = 2+2t^3+t^4 + \cdots + t^{n-1} \in \bZ[t]/(t^n-1)$$
in which the coefficient of $t^i$ is the number of times that
$z^i$ appears as a coordinate in~$x$. Then
$G_x(t)G_x(t^{-1}) \in \bZ[t]/(t^n-1)$ is the formal generating
function of $\rho(x)$.  Thus (ii) is equivalent to
\begin{equation} \label{e.gen-function}
\text{$G_x(t)G_x(t^{-1}) = G_y(t) G_y(t^{-1})$ in $\bZ[t]/(t^n -1)$.}
\end{equation}
Using the factorization
$t^n-1 = (t-1)(t^{n-1} +  \cdots +1)$, we get an isomorphism
\[ \text{$\bZ[t]/(t^n-1) \simeq \bZ \oplus R$, where
$R := \bZ[t]/(t^{n-1} + \cdots + 1)$.} \]
Here $t$ maps to~$1$ in the first summand.
We see that in order to prove~\eqref{e.gen-function}
it suffices to check that the two sides have the same images
in $\bZ$ and in $R$.
The image in $\bZ$ is obtained by evaluating at 1, and
$G_x(1) = G_y(1) = n$, so the desired identity is immediate in~$\bZ$.
If $g_x$ and $g_y$ denote the images of $G_x$
and $G_y$ in $R$, then a calculation similar to~\eqref{e.rt-of-unity}
show that $g_y(t) = -g_x(t)$.

It follows immediately that $g_x(t) g_x(t^{-1}) = g_y(t) g_y(t^{-1})$
as desired.  This finishes the proof for $n > 4$.

For $n = 4$ we have to consider a more carefully crafted
example.  Let $z$ be a primitive 13-th root of unity,
and set
\begin{equation} \label{e.diff-set}
x = (1 , z , z^4 , z^6), \qquad y = (1 , z , z^4 , z^{11}) \, .
\end{equation}
An easy calculation, which we will leave to the reader,
shows that $f(x) =  - f(y) \ne 0$, and a generating
function argument shows that $\rho(x) \sim \rho(y)$.
(This choice of $x$ and $y$ is based on the fact that
$\{ 0, 1, 4, 6 \}$ and
$\{ 0, 1, 4, 11 \}$ are inequivalent planar difference sets
for the cyclic group $\bZ/13 \bZ$; cf.~\cite[Chapter 9]{ryser}.)

This completes the proof of Proposition~\ref{prop.n>=5}.
\end{proof}

\section{Rationality}
\label{sect.rationality}

In this section we investigate the structure of the field
$k(X_1, \dots, X_m)_0^G$, without specific reference
to the observables. The main result, Proposition~\ref{prop.rationality}
below, was communicated to us by M. Lorenz.

Recall that a field extension $K/k$ is called {\em rational} (or
equivalently, $K$ is said to be rational over $k$)
if $K = k(t_1, \dots, t_r)$ for some
elements $t_1, \dots, t_r \in K$, algebraically independent over $k$.
The extension $K/k$ is called {\em stably rational} if
there exists a field $L$, containing $K$, which is rational over both
$K$ and $k$. In other words, for some finite collection of indeterminates
$s_1, \dots, s_r$ the field $L = K(s_1, \dots, s_r)$ is rational over $k$.

\begin{prop} \label{prop.rationality}
(a) For $m = 1$, $k(X)_0^G$ is rational over $k$.

\smallskip
(b) For any $m \geq 1$, $k(X_1, \dots, X_m)_0^G$ is stably
rational over $k$.
\end{prop}

\begin{proof} (a) First assume $n = 2$. Then
$k(X)^G \subset k(X)_0 = k(x_1/x_2)$, and the desired
conclusion follows from L\"uroth's theorem. For
$n \ge 3$, part (a) is a special case of a theorem of N. Lemire;
see~\cite[Theorem 7.7]{lemire}.  To see how
Lemire's theorem applies, note that the elements $x_i/x_j$
(viewed additively, with the natural $\Sym_n$-action)
form the root system $A_{n-1}$. The group $G = \Sym_n \times T$
is the automorphism group of this root system,
and $k(X)_0 = k(M)$, where $M$ is the lattice
(i.e., the multiplicative subgroup of $k(X)_0^*$) generated by
$\{ x_i/x_j \, | \, i, j = 1, \dots, n \}$.

\smallskip
(b) We argue by induction on $m$. The base case, $m = 1$, is part (a).
For the induction step, assume $m \geq 2$. Let $K_m = k(X_1, \dots, X_m)_0$
and $L = K_{m-1}(x_1, \dots, x_n)$, where $X_m = (x_1, \dots, x_n)$.
By the induction assumption, $K_{m-1}^G$ is stably rational over $k$;
thus it is enough to show that $K_m^G$ is stably rational
over $K_{m-1}^G$. We will do this by proving that $L^G$ is rational over
both $K_m^G$ and $K_{m-1}^G$.

Note that
since $L = K_m(x_1) = K_{m-1}(x_1, \dots, x_n)$ and the $G$-action on
$K_{m-1}$ is faithful, the desired conclusion
would follow from the no-name lemma, if $G$-acted linearly
(or, more precisely, semi-linearly) on the variables $x_1, \dots, x_n$;
see e.g.,~\cite[Proposition 1.1]{em} % or~\cite[Theorem 1]{hk1}
or~\cite[Appendix 3]{shafarevich}.
However, the no-name lemma cannot be used
directly in this setting because our $G$-action
is not semi-linear ($\tau$ acts by inversion!).

Fortunately, our action can be easily linearized, following an
argument of Hajja and Kang~\cite[Lemma 2.3(i)]{hk}.
Let $y_i = \frac{1-x_i}{1+x_i}$. Then $\Sym_n$
permutes $y_1, \dots, y_n$ and $\tau$ sends each $y_i$ to $-y_i$.
Since $x_i = \frac{1-y_i}{1+y_i}$, we have
\[ L = K_m(y_1) = K_{m-1}(y_1, \dots, y_n) \, . \]
Now the no-name lemma tells us that $L^G$ is rational over both
$K_m^G$ and $K_{m-1}^G$, as claimed.
\end{proof}

\begin{remark} \label{rem.hajja-kang} Let $H$ be the subgroup of
index 2 in $G = \Sym_n \times T$ consisting of pairs of the form
$(\sigma, \tau^{{\rm sign} \, {\sigma}})$. Hajja and
Kang~\cite[Theorem 3.2]{hk} have
shown that $k(X_1, \dots, X_m)_0^H$ is rational over $k$.
\end{remark}

\end{document}